\documentclass[reqno,a4paper]{amsart}
\usepackage{a4wide}
\usepackage[utf8x]{inputenc}
\usepackage{ucs}
\usepackage{amsmath}
\usepackage{amsthm}
\usepackage{amsfonts}
\usepackage{amssymb}
\usepackage{paralist}
\usepackage{graphicx}
\usepackage{xfrac}

\title[The dimension of the graph of the classical Weierstrass function]{An elementary proof for the dimension of the graph of the classical Weierstrass function}

\author{Gerhard Keller}
\address{Department Mathematik, Universit\"at Erlangen-N\"urnberg, 91058 Erlangen, Germany}
\email{keller@mi.uni-erlangen.de}
\thanks{This work is funded by DFG grant Ke 514/8-1. I am indebted to Tobias Oertel-J\"ager who not only awakened my interest in this problem but made it also a topic of the DFG Scientific Network ``Skew Product Dynamics and Multifractal Analysis'', and I want to thank Atsuya Otani for pointing out some mistakes in a previous version of this note.}
\date{\today}

\newtheorem{theorem}{Theorem}[section] % 1st argument is your name for it
\newtheorem{proposition}[theorem]{Proposition}

\theoremstyle{definition}

\newtheorem{remark}[theorem]{Remark}

\renewcommand{\mod}{\text{ mod }}
\newcommand{\A}{\mathbb{A}}
\newcommand{\E}{\mathbb{E}}
\renewcommand{\L}{\mathbb{L}}
\newcommand{\R}{\mathbb{R}}
\newcommand{\N}{\mathbb{N}}
\newcommand{\Z}{\mathbb{Z}}

\newcommand{\T}{{\mathbb{T}}}
\newcommand{\I}{{\mathbb{I}}}
\newcommand{\J}{{\mathbb{J}}}

\newcommand{\cG}{\mathcal{G}}

\newcommand{\cJ}{\mathcal{J}}

\newcommand{\x}{x}
\newcommand{\y}{y}

\newcommand{\const}{\text{const}}

\renewcommand{\u}{u}
\renewcommand{\v}{v}
\newcommand{\half}{\sfrac{1}{2}}
\newcommand{\btel}{\sfrac{1}{b}}
\numberwithin{equation}{section}

\begin{document}
\begin{abstract}
Let $W_{\lambda,b}(\x )=\sum_{n=0}^\infty\lambda^n g(b^n \x )$ where $b\geqslant2$ is an integer  and $g(u)=\cos(2\pi u)$ (classical Weierstrass function).
Building on work by Ledrappier \cite{Ledrappier1992}, Bar\'ansky, B\'ar\'any and Romanowska \cite{BBR2013}
and Tsujii \cite{Tsujii2001}, we provide an elementary proof that the Hausdorff dimension of $W_{\lambda,b}$ equals $2+\frac{\log\lambda}{\log b}$ for all $\lambda\in(\lambda_b,1)$ with a suitable $\lambda_b<1$. 
This reproduces results by 
Bar\'ansky, B\'ar\'any and Romanowska \cite{BBR2013} without using
the dimension theory for hyperbolic measures of Ledrappier and Young \cite{LY1985,LY1988}, which is replaced by a simple telescoping argument together with a recursive multi-scale estimate.
\end{abstract}

\maketitle

\section{Introduction}

The classical Weierstrass function $W_{\lambda,b}:\I:=[0,1)\to\R$ with parameters $b\in\N$, $\lambda\in(0,1)$ and $b\lambda>1$ is defined by
\begin{displaymath}
W_{\lambda,b}(\x )=\sum_{n=0}^\infty\lambda^n\cos(2\pi b^n \x )\ .
\end{displaymath}
The box dimension of its graph is equal to 
\begin{equation*}\label{eq:D-identity}
D=2+\frac{\log\lambda}{\log b}
\end{equation*}
as proved by Kaplan, Mallet-Paret and Yorke in \cite{KMY1984}. In 1977, Mandelbrot conjectured in his monograph \cite{Mandelbrot1977} that $D$ is also the Hausdorff dimension of this graph.
Despite many efforts, this conjecture is not yet proved in full generality.
Among others it is known to be true for sufficiently large integers \cite{Biacino2011,Fu2011}. The history of the problem and the present state of knowledge are summarized in the introduction to a recent paper by Bar\'ansky, B\'ar\'any and Romanowska \cite{BBR2013}, in which the authors
make a huge step forward by proving that for each integer $b\geqslant2$ there exist $\tilde{\lambda}_b<\lambda_b<1$ such that $D=2+\frac{\log\lambda}{\log b}$ equals the Hausdorff dimension of the graph of $W_{\lambda,b}$  for every $\lambda\in(\lambda_b,1)$ and for Lebesgue-a.e. $\lambda\in(\tilde{\lambda}_b,1)$. They determine $\lambda_b$ and $\tilde{\lambda}_b$ as unique zeroes of certain functions and provide a number of numerical and asymptotic values for them, among others
\begin{displaymath}
\lambda_2=0.9531,\quad\lambda_3=0.7269,\quad\lambda_4=0.6083,\quad\text{and}\quad\lim_{b\to\infty}\lambda_b=1/\pi=0.3183\ .
\end{displaymath}
Following Ledrappier \cite{Ledrappier1992}, they interpret the graph of $W_{\lambda,b}$ as the unique invariant repellor of the dynamical system
\begin{equation}\label{eq:weierstrass-repellor}
\Phi_{\lambda,b}:\I\times\R\to\I\times\R,\quad
\Phi_{\lambda,b}(\u,\v)=\left(b\u\mod1,\frac{\v-g(\u)}{\lambda}\right)
\end{equation}
with $g(\u)=\cos(2\pi\u)$, and observe that it suffices to show that $D$ is the Hausdorff dimension of the lift of the Lebesgue measure on $\I$ to the graph of $W_{\lambda,b}$, denoted by $\mu_{\lambda,b}$. Then
they  extend the transformation $\u\mapsto b\u\mod 1$ of the first coordinate to an invertible '$b$-baker' map.\footnote{More precisely, they consider the extension of $\u\mapsto b\u\mod 1$ by a full one-sided $b$-shift.} The resulting $3$-dimensional system is hyperbolic, and the extension of $\mu_{\lambda,b}$ is a hyperbolic invariant measure for it. This sets the stage to combine dimension results for hyperbolic measures by Ledrappier and Young \cite{LY1985} and \cite{LY1988}, an observation by Ledrappier \cite{Ledrappier1992} and a transversality estimate by Tsujii \cite{Tsujii2001} to determine $\lambda_b$. Finally, additional effort is needed to determine $\tilde{\lambda}_b$ based on the transversality approach of Peres and Solomyak \cite{PS1996}.

In this note I propose a much more elementary approach to reduce the calculation of the dimension of the graph of the Weierstrass function to the basic estimates provided by Tsujii \cite{Tsujii2001} combined with the numerical estimates by Bar\'ansky, B\'ar\'any and Romanowska \cite{BBR2013}. 
The proof also highlights conceptual similarities to the result by Hunt \cite{Hunt1998}, who studied randomized versions of the Weierstrass graph and proved the validity of the dimension formula for almost all random realizations, see Remark~\ref{remark:Hunt}.
An additional benefit of the present approach is that it avoids reference to an argument from \cite{Ledrappier1992}, that sketches only very briefly how to extend a 
certain lemma
for the piecewise linear function
$g(u)=\text{dist}(u,\Z)$ to nonlinear functions like $g(u)=\cos(2\pi u)$ by using
techniques from \cite{LY1988}. In passing we also consider the case of the piecewise linear function $g(u)=\text{dist}(u,\Z)$: we provide a short  and elementary argument how to reduce the calculation of the dimension of $W_{\lambda,b}$ in this case to the problem of whether an infinite Bernoulli convolution is absolutely continuous, a result that is originally proved in \cite{Ledrappier1992}. See Subsection~\ref{subsubsec:pwl} for details.

Although this note bypasses completely the general dimension theory of Ledrappier and Young, most proofs in this note are the result of my efforts to understand the basic lines of arguments in the papers mentioned above. The only, but noticeable, exception is Proposition~\ref{prop:inductive}, that  provides a new point of view on how to deal in a more direct way than in \cite{Ledrappier1992} with the passage from piecewise linear to nonlinear functions $g(u)$.

\section{The main results}
Throughout this note we use the notation
\begin{displaymath}
\gamma:=\frac{1}{b\lambda}<1\ .
\end{displaymath}
Our main results are the new proofs for the following two theorems - not the theorems themselves. The first one is due to Ledrappier \cite{Ledrappier1992}:

\begin{theorem}\label{theo:pwl}
Let $g(u)=\text{dist}(u,\Z)$, $b=2$, and let $\lambda\in(0,1)$ be such that the infinite Bernoulli convolution with parameter $\gamma$ has a square-integrable density w.r.t. Lebesgue measure. Then the graph of $W_{\lambda,2}$ 
has Hausdorff dimension $D=2+\frac{\log
\lambda}{\log 2}$.
\end{theorem}

\begin{remark}
The infinite Bernoulli convolution with parameter $\gamma$ is the distribution of the random variable $\Theta=\sum_{n=1}^\infty\gamma^nZ_n$, where the $Z_n$ are independent random variables with $P(Z_n=1)=P(Z_n=-1)=\frac{1}{2}$. The investigation of $\Theta$ has a long history, see e.g. \cite{Peres1999,Peres2000} and, for more recent results, also \cite{Shmerkin2013}. 

In particular, the set of parameters $\lambda\in(\half,1)$ for which the corresponding Bernoulli convolution with parameter $\gamma$ has a square integrable density, has full Lebesgue measure in this interval. (It corresponds to $\gamma\in(\half,1)$.) At the expense of only little additional effort our proof extends to the slightly more general case where  the distribution of $\Theta$ is only assumed to have dimension $1$, and also integers $b>2$ can be treated in just the same way.
\end{remark}

The second theorem is due to Barański, B\'{a}r\'{a}ny and Romanowska \cite{BBR2013}, building upon work of Ledrappier \cite{Ledrappier1992} and
Tsujii \cite{Tsujii2001}:

\begin{theorem}\label{theo:BBR}
Let $g(u)=\cos(2\pi u)$. For each integer $b\geqslant2$ there exists $\lambda_b<1$ such that the graph of $W_{\lambda,b}$ described by (\ref{eq:weierstrass-repellor}) has Hausdorff dimension $D=2+\frac{\log\lambda}{\log b}$  for every $\lambda\in(\lambda_b,1)$.
\end{theorem}

\begin{remark}
$\lambda_b$ is the unique zero of the function
\begin{displaymath}
h_b(\lambda)
=
\begin{cases}
\frac{1}{4\lambda^2(2\lambda-1)^2}+\frac{1}{16\lambda^2(4\lambda-1)^2}-\frac{5}{64\lambda^2}+\frac{\sqrt 2}{2\lambda}-1&\text{for }b=2\\
\frac{1}{(b\lambda-1)^2}+\frac{1}{(b^2\lambda-1)^2}-\sin^2(\frac{\pi}{b})
&\text{for }b\geqslant3
\end{cases}
\end{displaymath}
on the interval $(\btel,1)$, see \cite[Theorem B]{BBR2013}.
\end{remark}

\section{Proofs}
In Sections~\ref{subsec:Weierstrass-attractor} and~\ref{subsec:stable-unstable} we recall some observations from \cite{Ledrappier1992} and \cite{BBR2013}, and in Section~\ref{subsec:distances} we provide a fresh look at the strong stable manifolds from those references. Section~\ref{subsec:telescoping} contains the telescoping argument already used in a
similar situation in \cite{Keller2012a},
 and the proof is finished in Sections~\ref{subsec:Marstrand} and~\ref{subsec:last-proof} by combining some of the more elementary arguments from \cite{Ledrappier1992}, \cite{Tsujii2001} and \cite{BBR2013}.
\subsection{The Weierstrass graph as an attractor}\label{subsec:Weierstrass-attractor}
Recall from (\ref{eq:weierstrass-repellor}) that 
\begin{displaymath}
\Phi_{\lambda,b}(\u,\v)=\left(b\u\mod1,\frac{\v-g(\u)}{\lambda}\right)\ .
\end{displaymath}
We are mostly interested in the classical case $g(u)=\cos(2\pi u)$ and in 
$g(u)=\text{dist}(u,\Z)$
where $g'(\xi )=(-1)^{\lfloor 2\xi \rfloor}$.
For notational convenience we denote the map $\u\mapsto b\u\mod 1$ by $\tau$ so that
$\Phi_{\lambda,b}(\u,\v)=(\tau(\u),\frac{\v-g(\u)}{\lambda})$.
 Then the Weierstrass function $W=W_{\lambda,b}$ satisfies
\begin{displaymath}
\Phi(\u,W(\u))=(\tau(\u),W(\tau(\u)))\ .
\end{displaymath}
In particular, 
\begin{equation}\label{eq:weierstrass-2}
\lambda W(\tau(\u))=W(\u)-g(\u)\ .
\end{equation}

Denote by $(\xi ,\x )\mapsto B(\xi ,\x )$ the $b$-baker map on $\I^2$ for the integer $b\geqslant 2$, i.e.
\begin{displaymath}
B(\xi ,\x )=\left(\tau(\xi ),\frac{\x +k(\xi )}{b}\right)\quad\text{with}\quad
k(\xi )=j\in\{0,\dots,b-1\}\text{ if }\xi \in[j/b,(j+1)/b)\ ,
\end{displaymath}
and define $F:\I^2\times\R\to\I^2\times\R$ as
\begin{equation*}
F(\xi ,\x ,\y )=\left(B(\xi ,\x ),\lambda \y +f(\xi ,\x )\right)\quad\text{with}\quad
f(\xi ,\x ):=g\left(\frac{\x +k(\xi )}{b}\right)\ .
\end{equation*}
Then the graph of the Weierstrass function $W=W_{\lambda,b}$, interpreted as a function of $\xi$ and $x$ although it depends on $x$ only, is an invariant attractor for $F$ in the following sense:
\begin{displaymath}
\begin{split}
F(\xi ,\x ,W(\xi,x ))
&=
\left(B(\xi,x),\lambda W(\x )+g\left(\frac{\x +k(\xi )}{b}\right)\right)
=
\left(B(\xi,x), W\left(\frac{\x +k(\xi )}{b}\right)\right)\\
&=
\left(B(\xi ,\x ), W\left(B(\xi ,\x )\right)\right)\ ,
\end{split}
\end{displaymath}
where the second identity follow from (\ref{eq:weierstrass-2}) with $\u=\frac{\x +k(\xi )}{b}$. As $F$ has skew-product structure over the base $B$ and as $\left|\frac{\partial F_3}{\partial \y }\right|=\lambda<1$, the graph of $W$ is an attractor for $F$.

\subsubsection{Notation for orbits}
Given a point $(\xi ,\x )\in\I^2$ we denote by $(\xi _n,\x _n)$ the point $B^n(\xi ,\x )$ ($n\in\Z$).
Note that $\xi _n=\tau^n(\xi )$ ($n\geqslant0$) and $\x _n=\tau^{-n}(\x) $ ($n\leqslant0$) and that
\begin{equation*}
\label{eq:x_i-y_i}
k(\xi _i)=k(\x _{i+1})\quad\text{for all }i\in\Z\ .
\end{equation*}
We also write
\begin{displaymath}
k_n(\xi )=\sum_{i=0}^{n-1} b^{i} k(\tau^{i}\xi )\ .
\end{displaymath}
For later use we note that 
\begin{displaymath}
\begin{split}
k_n(\xi )
&=
\sum_{i=0}^{n-1}b^ik(\xi _i)
=
b^n\,\sum_{i=0}^{n-1}b^{i-n}k(\x _{i+1})
=
b^n\,\sum_{j=0}^{n-1}b^{-j-1}k(\x _{n-j})
=
b^n\,\sum_{j=0}^{n-1}b^{-j-1}k(\tau^j\x _{n})\\
&=
b^n\,\x _n-b^n\sum_{j=n}^\infty b^{-j-1}k(\x _{n-j})
=
b^n\,\x _n-\sum_{i=0}^\infty b^{-i-1}k(\x _{-i})
=
b^n\,\x _n-\sum_{i=0}^\infty b^{-i-1}k(\tau^i\x )\\
&=
b^n\,\x _n-\x \ .
\end{split}
\end{displaymath}
In particular,
\begin{equation*}
\label{eq:knx}
\begin{split}
k_n(\xi _{-i})&=b^n\,\tau^{i-n}(\x )-\tau^i(\x )\quad\text{for }n\leqslant i\ ,\\
k_n(\xi _{-i})&=b^n\,\x _{n-i}-\tau^i(\x )\hspace{8.4mm}\text{for }n\geqslant i\ ,
\end{split}
\end{equation*}
and
\begin{equation*}
\frac{\x_n+k(\xi_n)}{b}
=
\frac{\x+k_n(\xi)+b^{n}\,k(\xi_n)}{b^{n+1}}
=
\frac{\x+k_{n+1}(\xi)}{b^{n+1}}
=
x_{n+1}\ .
\end{equation*}
For comparison with the notation of \cite{BBR2013} note also that
\begin{equation}
\label{eq:xk4}
\x_n
=
\frac{x+k_n(\xi)}{b^n}
=
\frac{\x}{b^n}+\frac{k(\x_1)}{b^n}+\dots+\frac{k(\x_n)}{b}
=
\frac{\x}{b^n}+\frac{k(\xi_0)}{b^n}+\dots+\frac{k(\xi_{n-1})}{b}\ .
\end{equation}

\subsection{Stable and unstable manifolds}\label{subsec:stable-unstable}
Following \cite{BBR2013} and also the earlier paper \cite{Ledrappier1992}, we describe the stable and unstable manifolds of $F$. The derivative $DF$ is well defined except when $\xi \in S:=\{j/b: j=0,\dots,b-1\}$, namely
\begin{displaymath}
DF(\xi ,\x ,\y )
=
\begin{pmatrix}
b&0&0\\
0&\frac{1}{b}&0\\
0&\frac{\partial f}{\partial \x }(\xi ,\x )&\lambda
\end{pmatrix}
.
\end{displaymath}
The Lyapunov exponents of the corresponding cocycle are $\log b$, $-\log b$ and $\log\lambda$. Indeed, they correspond to the invariant vector fields
\begin{displaymath}
\begin{pmatrix}
1\\
0\\
0
\end{pmatrix},\quad
X(\xi ,\x ,\y )
=
\begin{pmatrix}
0\\1\\
-\lambda^{-1}\sum_{n=0}^\infty\gamma^n\frac{\partial f}{\partial \x }(B^{n}(\xi ,\x ))
\end{pmatrix},\text{ and }\quad
\begin{pmatrix}
0\\0\\1
\end{pmatrix}
\end{displaymath}
where $\gamma=(b\lambda)^{-1}$.
Observe that none of these fields depends on the variable $\y $, so we write
$X(\xi ,\x )$ henceforth. 
\begin{remark}
As $\frac{\partial f}{\partial \x }(B^{n}(\xi ,\x ))=b^{-1}g'\left(\frac{\x_n+k(\xi_n)}{b}\right)
=b^{-1}g'\left(x_{n+1}\right)$, the third component of $X$ can be written as
\begin{displaymath}
X_3(\xi,x)
=
-\sum_{n=1}^\infty\gamma^n\,g'\left(x_n\right)
=
-\sum_{n=1}^\infty\gamma^n\,g'\left(\frac{\x+k_n(\xi)}{b^n}\right)
%=
%-\sum_{n=1}^\infty\gamma^n\,g'\left(\frac{\x}{b^n}+\frac{k(\x_1)}{b^n}+\dots+\frac{k(\x_n)}{b}\right)
\end{displaymath}
which is precisely the second component of the field $\cJ_{x,{\bf i}}$ of \cite{BBR2013}, also denoted $Y_{x,\gamma}({\bf i})$ in that paper.
\end{remark}

For each fixed $\xi $, the field $X$ defines the strong stable foliation in the $(\x ,\y )$-plane $H_\xi $ over $\xi $. The fibres are parallel graphs over $\x $ with uniformly bounded slopes. 
That means, for all $(\xi ,\x ,\y )\in(\I\setminus S)\times\I\times\R$ the fibre through $(\xi ,\x ,\y )$ is the graph of a function $\ell^{ss}_{(\xi ,\x ,\y )}:\I\to\R$ that solves the initial value problem
\begin{displaymath}
\frac{\partial}{\partial v}\ell^{ss}_{(\xi ,\x ,\y )}(v)=X_3(\xi ,v)
\quad\text{and}\quad
\ell^{ss}_{(\xi ,\x ,\y )}(\x )=\y \ .
\end{displaymath}

Denote by $\cG\ell^{ss}_{(\xi ,\x ,\y )}$ the graph of the function $\ell^{ss}_{(\xi ,\x ,\y )}$ in the hyperplane $H_\xi $, i.e.
\begin{equation*}
\cG\ell^{ss}_{(\xi ,\x ,\y )}=\left\{(\xi ,\u,\ell^{ss}_{(\xi ,\x ,\y )}(\u)):\,\u\in\I\right\}\ .
\end{equation*}
As the foliation into strong stable fibres is invariant, we have
\begin{equation*}
F\left(\cG\ell^{ss}_{(\xi ,\x ,\y )}\right)
\subseteq
\cG\ell^{ss}_{F(\xi ,\x ,\y )}\ .
\end{equation*}

\subsection{Distances between strong stable fibres}\label{subsec:distances}
Given two points $(\xi ,\x ),(\xi ,\x ')\in\I^2$ we denote by $|\Delta_\xi (\x ,\x ')|$ the vertical distance of the strong stable fibres through the points $(\xi ,\x ,W(\x ))$ and $(\xi ,\x ',W(\x '))$, respectively.
More precisely,
\begin{equation*}
\begin{split}
\Delta_\xi (\x ,\x ')
&=
\ell^{ss}_{(\xi ,\x ',W(\x '))}-\ell^{ss}_{(\xi ,\x ,W(\x ))}
=
\ell^{ss}_{(\xi ,\x ',W(\x '))}(x')-\ell^{ss}_{(\xi ,\x ,W(\x ))}(x')
\\
&=
W(\x ')-W(\x )-\left(\ell^{ss}_{(\xi ,\x ,W(\x ))}(\x ')-\ell^{ss}_{(\xi ,\x ,W(\x ))}(\x )\right)\\
&=
W(\x ')-W(\x )-\int_\x ^{\x '}X_3(\xi ,t)\,dt\\
&=
W(\x ')-W(\x )+\sum_{n=1}^\infty\gamma^n\int_x^{x'}g'\left(\frac{t+k_n(\xi)}{b^n}\right)dt
\end{split}
\end{equation*}

\subsubsection{The piecewise linear case}
\label{subsubsec:pwl}
In the piecewise linear case
$g'(u)=(-1)^{\lfloor 2u \rfloor}$ and $b=2$ we have
\begin{displaymath}
g'\left(\frac{t+k_n(\xi)}{b^n}\right)
=
g'\left(\frac{t}{b^n}+\frac{k(\xi_0)}{b^n}+\dots+\frac{k(\xi_{n-1})}{b}\right)
=
g'\left(\frac{k(\xi_{n-1})}{2}\right)=(-1)^{k(\xi_{n-1})}
\end{displaymath}
so that
\begin{equation}
\begin{split}
\Delta_\xi(\x,\x')
&=
W(\x')-W(\x)+
(\x'-\x)\cdot\sum_{n=1}^\infty\gamma^n\,(-1)^{k(\xi _{n-1})}\\
&=
W(\x')-W(\x)+
(\x'-\x)\cdot\Theta(\xi )
\end{split}
\end{equation}
where $\Theta(\xi ):=\sum_{n=1}^\infty\gamma^n\,(-1)^{k(\xi _{n-1})}$ is an infinite Bernoulli convolution. It is known \cite{Peres1999,Shmerkin2013} that 
for Lebesgue-a.e. $\gamma\in(\half,1)$ the distribution of the random variable $\Theta$ has a square-integrable density $h$ w.r.t. Lebesgue measure. For such $\Theta$ the following holds:

 For each $\delta>0$ there is $C>0$ such that for each $r>0$ there is a measurable set $E_r\subset\I^2$ with $m^2(E_r)\leqslant C r^\delta$ and the following property: For each $z\in[-1,1]\setminus\{0\}$ and for each
measurable family
$(\L_x)_{x\in\I}$ of intervals of length $\frac{2r}{|z|}$,
\begin{equation}\label{eq:main-prop-pwl}
m^2\{(\xi,x)\in \I^2\setminus E_r: \Theta_z(\xi,x)\in \L_x\}
\leqslant
r^{1-2\delta}\,|z|^{-1}
\ .
\end{equation}
Indeed, let $E_r=\{(\xi,x)\in\I^2: h(\Theta(\xi,x))>r^{-\delta}\}$.
Then $m^2(E_r)=\int_\I\int_{\{\theta:h(\theta)>
r^{-\delta}\}}h(\theta)\,d\theta\,dx\leqslant r^{\delta}\|h\|_2^2$
and
\begin{displaymath}
m^2\{(\xi,x)\in \I^2\setminus E_r: \Theta_z(\xi,x)\in \L_x\}
=
\int_{\{\theta:h(\theta)\leqslant r^\delta\}}1_{\L_x}(\theta)\,h(\theta)\,d\theta
\leqslant
r^{1-\delta}|z|^{-1}\ .
\end{displaymath}
This is a simplified version of the much more difficult Proposition~\ref{prop:inductive} that we state and prove below for the case of nonlinear $g$. Readers who want first to see how the proof for the present case is finished, can jump immediately to Section~\ref{subsec:telescoping} and only read that one and Section~\ref{subsec:Marstrand}.

\subsubsection{The case $g(u)=\cos(2\pi u)$}
If $g(u)=\cos(2\pi u)$, then
\begin{displaymath}
\begin{split}
\int_x^{x'}g'\left(\frac{t+k_n(\xi)}{b^n}\right)dt
&=
b^n\left(\cos\left(2\pi\frac{x'+k_n(\xi)}{b^n}\right)-\cos\left(2\pi\frac{x+k_n(\xi)}{b^n}\right)\right)\\
&=
-2b^n\sin\left(2\pi\frac{x'-x}{2b^n}\right)\sin\left(2\pi\left(\frac{x'+x}{2}+
k_n(\xi)\right)\big/b^n\right)
%&=
%-2b^n\sin\left(\pi\frac{x'-x}{b^n}\right)\sin\left(2\pi\left(x+\frac{x'-x}{2}+
%k_n(\xi)\right)\big/b^n\right)
\ ,
\end{split}
\end{displaymath}
so that, with
$s(t):=(t/2)^{-1}\sin(2\pi\, t/2)$,
\begin{displaymath}
\Delta_\xi (\x ,\x ')
=
W(x')-W(x)-(x'-x)\sum_{n=1}^\infty\gamma^n s\left(\frac{x'-x}{b^n}\right)
\sin\left(2\pi\left(x+\frac{x'-x}{2}+
k_n(\xi)\right)\big/b^n\right)\ .
\end{displaymath}
With
\begin{equation*}
\Theta_z(\xi,x):=\sum_{n=1}^\infty\gamma^n\,s\left(\frac{z}{b^n}\right)
\sin\left(2\pi\left( x_n+\frac{z}{2b^n}\right)\right)
\end{equation*} 
this can be written as
\begin{equation}
\Delta_\xi (\x ,\x ')
=
W(x')-W(x)-(x'-x)\cdot\Theta_{x'-x}(\xi,x)\ ,
\end{equation}
because $\frac{x+k_n(\xi)}{b^n}=x_n$, see (\ref{eq:xk4}).

The function $\Theta_0(\xi,x)=-2\pi\sum_{n=1}^\infty\gamma^n\sin(2\pi x_n)$ is, up to some constant factor and different notation, just the function $S(x,\mathbf{i})$ from \cite{BBR2013}, and Proposition 4.2 of this paper (which is proved via some explicit estimates) together with the more elementary part of Tsujii's paper \cite[Sections 3, 4]{Tsujii2001} yields the following fact:
Denote by
$\nu_{x}:=(m\times\delta_x)\circ\Theta_0^{-1}$
the conditional distribution of $\Theta_0(\xi,x)$ given $x\in\I$.
\begin{proposition}\label{prop:like-Tsujii}
Let $\lambda\in(\lambda_b,1)$, i.e. $b\gamma\in(1,\lambda_b^{-1})$.
For $m$-a.e. $x\in\I$, the measure $\nu_{x}$ is absolutely continuous w.r.t. $m$. Its density $h_{x}$ satisfies
$H:=\int_{\I}\|h_{x}\|_2^2\,dx<\infty$.
\end{proposition}
A major technical problem is that this estimate is needed also for $z\neq0$. One approach could be to imitate Tsujii's recursion from \cite{Tsujii2001}, and indeed, one obtains densities $h_{x}$ with $\sup_{|z|\leqslant1}\int_\I\|h_{x}\|_2^2\,dx\\<\infty$. But this approach does not provide any local information on the $h_{x}$ uniformly in $z$: the set of $(\xi,x)$ where $h_{x}(\xi)$ is exceptionally big, might depend on $z$ in a complicated way. Therefore we follow a different approach here. Naively, one can start with comparing $\Theta_z$ to $\Theta_0$: it is easily seen that there is a constant $C>0$ such that $\|\Theta_z-\Theta_0\|_\infty\leqslant C|z|$. As in later steps of the proof we have to approximate $\Theta_z$ by $\Theta_0$ up to an error of order $r$ for small $r>0$, this would cover only $|z|<{r}$. \footnote{The same problem occurs also in Ledrappier's sketch of a related proof \cite{Ledrappier1992}. He solves it by using formulas relating dimensions and exponents of various conditional and projected measures as in \cite{LY1985}.}
 However, if one treats a finite part of the sum defining $\Theta_z$ separately from the remaining tail, one sees that a tail starting at $n=n_0$ varies with $z$ only of the order $\left(\frac{\gamma}{b}\right)^{n_0}|z|$. Using this observation recursively we will prove the following result in Section~\ref{subsec:last-proof}.

\begin{proposition}\label{prop:inductive}
Let $\lambda\in(\lambda_b,1)$. For each $\eta>0$ there are $\delta\in(0,\eta)$ and $C>0$ such that for each $r\in(0,1)$ there is a measurable set $E_r\subset\I^2$ with $m^2(E_r)\leqslant C r^\delta$ and the following property: For each $z\in[-1,1]\setminus\{0\}$ and for each
measurable family
$(\L_x)_{x\in\I}$ of intervals of length $\frac{2r}{|z|}$,
\begin{equation}\label{eq:main-prop}
m^2\{(\xi,x)\in \I^2\setminus E_r: \Theta_z(\xi,x)\in \L_x\}
\leqslant
C\, r^{1-2\eta}\,|z|^{-(1-\eta)}
\ .
\end{equation}
\end{proposition}
A crucial ingredient of the proof is the following observation:
\begin{remark}\label{remark:conditionalTheta}
Recall from (\ref{eq:xk4}) that $x_n=\frac{\x}{b^n}+\frac{k(\xi_0)}{b^n}+\dots+\frac{k(\xi_{n-1})}{b}$. Hence the conditional distribution,
 given $(x,k(\xi_0),\dots,k(\xi_{N-1}))$,
of 
\begin{equation*}
\begin{split}
\Theta_0(B^N(\xi,x))
&=
\Theta_0(\xi_N,x_N)=
-2\pi\sum_{n=1}^\infty\gamma^n\sin(2\pi x_{N+n})\\
&=
-2\pi\sum_{n=1}^\infty\gamma^n\sin\left(2\pi\left(\frac{x_N}{b^n}+\frac{k(\xi_N)}{b^n}+\frac{k(\xi_{N+1})}{b^{n-1}}+\dots+\frac{k(\xi_{N+n-1})}{b}\right)\right)
\end{split}
\end{equation*}
is $\nu_{x_N}$, the distribution of $\Theta_0(\,.\,,x_N)$, because the $k(\xi_n)$ are independent and uniformly distributed on $\{0,\dots,b-1\}$.
\end{remark}

\subsection{Telescoping - a replacement for the Ledrappier-Young argument}
\label{subsec:telescoping}

\subsubsection{Neighbourhoods bounded by strong stable fibres}
We define a kind of $\epsilon$-neighbourhoods of points $(\xi ,\x ,W(\x ))$ in $(\x ,\y )$-direction. To that end fix a constant $K>0$ (to be determined later) and, for any $\xi \in\I$ and a $b$-adic $\epsilon$-neighbourhood $I_N(\x )$ of $\x \in\I$ with $\epsilon=b^{-N}$, let
\begin{displaymath}
\begin{split}
V_N(\xi ,\x )
&=\left\{
(v,w)\in\I\times\R:\ v\in I_{N}(\x ),
|w-\ell^{ss}_{(\xi ,\x ,W(\x ))}(v)|\leqslant Kb^{-N}
\right\}\ .
\end{split}
\end{displaymath}
The sets $\{\xi \}\times V_N(\xi ,\x )$ are quadrilaterals in $H_\xi $, which are bounded in $\x $-direction by two vertical lines of distance $b^{-N}$ and in $\y $-direction by the strong stable fibres through
$(\xi ,\x ,W(\x )\pm Kb^{-N})$ (which are parallel!).
Denote by $G:=\{(\x ,W(\x )):\ \x \in\I\}$ the graph of $W$, let
\begin{displaymath}
A_N(\xi ,\x )=V_N(\xi ,\x )\cap G
\end{displaymath}
and let $\mu$ be the Lebesgue measure $m$ on $\I$ lifted to $G$. We will evaluate the local dimension (in $H_\xi $) of $\mu$ at $(\x , W(\x ))\in G$ along $b$-adic neighbourhoods $V_N(\xi ,\x )$,
i.e. we are going to determine the limit
\begin{equation}\label{eq:dim-limit-1}
\lim_{N\to\infty}\frac{\log\mu(V_N(\xi ,\x ))}{\log(b^{-N})}\ .
\end{equation}
Observe that this limit, if it exists, does not depend on $\xi $, as the next remark shows among others.
\begin{remark}
As $X_3$ is uniformly bounded by some constant $K_1$, all $\ell^{ss}_{(\xi ,\x ,W(\x ))}$ have $K_1$ as a common Lipschitz constant. Fixing the constant $K$ as $K_1+1$ and choosing $n_1\in\N$ such that $b^{n_1}>2K_1+1$, elementary geometric arguments show that
\begin{displaymath}
V_{N+n_1}(\xi ,\x )
\subseteq
\left\{(v,w)\in\I\times\R:\ v\in I_N(x), |w-W(x)|\leqslant b^{-N}\right\}
\subseteq
V_N(\xi,x)\ .
\end{displaymath}
This proves not only that the limit in (\ref{eq:dim-limit-1}) does not depend on $\xi$, but also that the $V_N(\xi,x)$ can be replaced by rectangles of height $2\cdot2^{-N}$ over the base $I_N(x)$.

Furthermore, for $m$-a.e. $x$,  one can replace
%It is well known that, the limit in (\ref{eq:dim-limit-1}) does not change if 
the dyadic intervals $I_N(x)$ by symmetric intervals
$I'_N(x):=[x-2^{-N},x+2^{-N}]$ and hence $V_N(\xi,x)$ by $V_N'(\xi,x):=I_N'(x)\times I_N'(W(x))$. Indeed, it is immediate that $V_{N+n_1}(\xi,x)\subseteq V'_N(\xi,x)$ and, by Borel-Cantelli, for $m$-a.e. $x$ there is $N(x)\in\N$ such that $V'_{N+[2\log_2 N]}(\xi,x)\subseteq V_N(\xi,x)$ for all $N\geqslant N(x)$.
\end{remark}

\subsubsection{The telescoping step}
$F^{-N}(\{\xi \}\times V_N(\xi ,\x ))$ is the image of the quadrilateral
$\{\xi \}\times V_N(\xi ,\x )$ in 
$H_{\xi _{-N}}$ under a map with derivative $\text{diag}(b^N,\lambda^{-N})$ which maps strong stable fibres to strong stable fibres. Hence
\begin{displaymath}
F^{-N}(\{\xi \}\times V_N(\xi ,\x ))
=
\{\xi _{-N}\}\times \Sigma_N(\xi _{-N},\x _{-N})
\end{displaymath}
where
\begin{displaymath}
\begin{split}
&\Sigma_N(\xi,\u)
:=
\left\{(v,w)\in\I\times\R:\ 
|w-\ell^{ss}_{(\xi,\u,W(\u))}(v)|\leqslant K(b\lambda)^{-N}
\right\}
\end{split}
\end{displaymath}
is a strip in $H_{\xi}$ of width $1$ and height $2K(b\lambda)^{-N}=2K\gamma^N$. Therefore,
\begin{equation}
%\label{eq:muVN}
\begin{split}
\frac{\mu(V_N(\xi ,\x ))}{m(I_N(\x))}
&=
\frac{m\left(\left\{v\in I_{N}(\x ):\,
|W(v)-\ell^{ss}_{(\xi ,\x ,W(\x ))}(v)|\leqslant Kb^{-N}
\right\}\right)}{m(I_N(\x))}\\
&=
\frac{m\left(\left\{\x'\in\I:\,|W(\x')-\ell^{ss}_{(\xi_{-N},\x_{-N},W(\x_{-N}))}(\x')|\leqslant K\gamma^{N}\right\}\right)}{m(\I)}\\
&=
\frac{m\left(\left\{\x'\in\I:\,|\Delta_{\xi _{-N}}(\x _{-N},\x ')|\leqslant K\gamma^{N}\right\}\right)}{m(\I)}
\end{split}
\end{equation}
so that
\begin{equation}
%\label{eq:lim-identity}
\lim_{N\to\infty}\frac{\log\mu(V_N(\xi,\x))}{\log(b^{-N})}
=
1+\lim_{N\to\infty}\frac{\log m\left\{\x'\in\I:\ |\Delta_{\xi_{-N}}(\x_{-N},\x')|\leqslant K\gamma^{N}\right\}}{\log(b^{-N})}
\end{equation}
provided the limits exist.
This corresponds to identity (2.3) in \cite{BBR2013}, which states 
that $\operatorname{dim}\mu=1+\frac{\log\gamma}{\log b^{-1}}\cdot\operatorname{dim}\nu_{x,{\bf i}}$ for certain probabilities $\nu_{x,{\bf i}}$.
Indeed, the remaining task in that paper, namely to show that $\operatorname{dim}\nu_{x,{\bf i}}\geqslant1$, corresponds in our approach to showing
that
\begin{equation}
\label{eq:intermediate}
\liminf_{N\to\infty}\frac{\log m\left\{\x'\in\I:\ |\Delta_{\xi_{-N}}(\x_{-N},\x')|\leqslant K\gamma^{N}\right\}}{\log(\gamma^{N})}
\geqslant 1\quad\text{for $m^2$-a.e. $(\xi,x)\in\I^2$.}
\end{equation}
We prove this in Section~\ref{subsec:Marstrand}.
%It can be interpreted in the following way:
%for ''typical'' $(\xi,\x)$ the distribution of the random variable
%$\Delta_{\xi_{-N}}(\x_{-N},\,.\,)$ has local dimension (at least) $1$ at $\Delta=0$. 
Indeed, for typical $(\xi,x)$, the distribution of the random variable
$\Delta_{\xi_{-N}}(\x_{-N},\,.\,)$ is closely related to the $\nu_{x,{\bf i}}$ of \cite{BBR2013}.

\begin{remark}\label{remark:Hunt}
Instead of projecting along stable fibres $\ell^s_{(\xi,x,W(x))}$ that depend on the additional variable $\xi$, one could as well choose a new coordinate system for each $\xi$, describe the Weierstrass function $W$ in this new coordinate system (resulting in a transformed version $W_\xi$ of $W$) and project the 
lifts of
Lebesgue measure to the graphs of the $W_\xi$, horizontally to the real axis. 
These projected measures would typically be different one from each other (they depend on $\xi$), but the arguments above show that they all have the same dimension. In this sense our approach is equivalent to determining the dimension of the graph of $W_\xi$ for almost all realisations of this random collection of graphs. For Weierstrass graphs with random phase shifts this was done by Hunt \cite{Hunt1998}. The difference to our situation is that Hunt introduced additional external randomness to the problem so that, for exceptional realizations, his random graphs may have a dimension different from the typical one. In contrast to that, in our case the randomness is generated by the dynamics itself, namely by the unstable coordinate of the underlying baker map,  and the construction guarantees that there are no exceptional realizations.
\end{remark}

\subsection{A Marstrand projection estimate}\label{subsec:Marstrand}
For our further discussion we use the assumption, covering both theorems, that the parameter $\gamma$ is such that 
the random variables $\Theta_z$ on $(\I^2,m^2)$ have  distributions of dimension $1$ in the sense that they obey the conclusion of Proposition~\ref{prop:inductive}. For the classical Weierstrass function with $\lambda \in(\lambda_b,1)$ we prove this proposition at the end of this note.
In the piecewise linear case where $\Theta_z(\xi,x)=\Theta(\xi)$ is an infinite Bernoulli convolution, this is an additional assumption satisfied for Lebesgue-almost $\gamma\in(\half,1)$ as discussed around Equation (\ref{eq:main-prop-pwl}).

The following argument is inspired by \cite{Ledrappier1992}. Let
$\eta>0$ and let $\delta\in(0,\eta)$, $C>0$ and the sets $E_r\subseteq\I^2$ be as in Proposition~\ref{prop:inductive}. 
Let $\A=\{(\xi,x,z)\in\I^2\times[-1,1]: 0\leqslant x+z\leqslant1\}$ and
$J_{r,x,z}:=\left[\frac{W(x+z)-W(x)}{z}-\frac{r}{|z|},\frac{W(x+z)-W(x)}{z}+\frac{r}{|z|}\right]$. Then
\begin{equation*}
\begin{split}
&\hspace{-0.5cm}
m^3\{(\xi ,\x ,\x ')\in\I^3: (\xi,x)\not\in\E_r, 
|\Delta_\xi (\x ,\x ')|\leqslant r\}\\
&=
m^3\left\{(\xi ,\x ,\x ')\in\I^3: (\xi,x)\not\in\E_r, 
|\Theta_{x'-x}(\xi,x)\cdot(\x '-\x )-(W(\x')-W(\x ))|\leqslant r\right\}\\
&=
m^3\left\{(\xi ,\x ,z)\in\A: (\xi,x)\not\in\E_r,  
\Theta_z(\xi,x)\in J_{r,x,z}\right\}\\
&=
\int_{-1}^1 m^2\left\{(\xi,x)\in\I^2\setminus E_r: (\xi,x,z)\in\A, \Theta_z(\xi,x)\in J_{r,x,z}\right\}\,dz\\
&\leqslant
C\int_{-1}^1r^{1-2\eta}|z|^{-(1-\eta)}\,dz
\hspace*{2cm}\text{(by Proposition~\ref{prop:inductive})}
\\
&\leqslant
C r^{1-2\eta}\ .
\end{split}
\end{equation*}
Here and in the sequel,
$C$ denotes a generic constant whose value may change from occurrence to occurrence and depend on the fixed quatities $\eta,\delta$ and $\ell$.
Therefore, writing again $(\xi_{-N},x_{-N})$ for $B^{-N}(\xi,x)$ and using the $B$-invariance of $m^2$,
\begin{displaymath}
\begin{split}
&\hspace*{-0.5cm}
m^2\left\{(\xi,x)\in\I^2: m\{x'\in\I: |\Delta_{\xi_{-N}}(x_{-N},x')|\leqslant r\}\geqslant r^{1-3\eta}\right\}\\
&=
m^2\left\{(\xi,x)\in\I^2: m\{x'\in\I: |\Delta_{\xi}(x,x')|\leqslant r\}\geqslant r^{1-3\eta}\right\}\\
&\leqslant
m^2(E_r)+
m^2\left\{(\xi,x)\in\I^2\setminus E_r: m\{x'\in\I: |\Delta_{\xi}(x,x')|\leqslant r\}\geqslant r^{1-3\eta}\right\}\\
&\leqslant
Cr^\delta+
r^{-(1-3\eta)}\int_{\I^2\setminus E_r}m\{x'\in\I: |\Delta_{\xi}(x,x')|\leqslant r\}\,dm^2(\xi,x)
\\
&=
Cr^\delta+
r^{-(1-3\eta)}\,m^3\{(\xi ,\x ,\x ')\in\I^3:\,(\xi,x)\not\in E_r, |\Delta_\xi (\x ,\x ')|\leqslant r\}\\
&\leqslant
Cr^\delta+Cr^{-(1-3\eta)}r^{1-2\eta}\\
&\leqslant
Cr^\delta\ .
\end{split}
\end{displaymath}
By Borel-Cantelli we thus conclude with $r=K\gamma^N$ that
\begin{displaymath}
\limsup_{N\to\infty}\gamma^{-(1-3\eta)N}
m\{x'\in\I: |\Delta_{\xi_{-N}}(x_{-N},x')|\leqslant K\gamma^N\}
\leqslant K^{1-3\eta}
\end{displaymath}
for $m^2$-a.e. $(\xi,x)\in\T^2$.
On a logarithmic scale this implies
\begin{displaymath}
\liminf_{N\to\infty}\frac{\log m\left\{x'\in\T^1:\ |\Delta_{\xi_{-N}}(x_{-N},x')|\leqslant K\gamma^{N}\right\}}{\log(b^{-N})}
\geqslant (1-3\eta)\,\frac{\log\gamma}{\log b^{-1}}\ ,
\end{displaymath}
and as this holds for all $\eta>0$, it proves (\ref{eq:intermediate}) and thus finishes the proofs of Theorems~\ref{theo:pwl} and~\ref{theo:BBR}.

\subsection{Proof of Proposition~\ref{prop:inductive}}
\label{subsec:last-proof}
Let $\eta>0$ and $\alpha:=\frac{\log\gamma}{\log(\gamma/ b)}<1$. Choose $\ell\in\N$ such that $\alpha^\ell<\eta$ and let $\delta=\frac{\eta}{2\ell}$. 

Given $r\in(0,1)$ and $z\in[-1,1]\setminus\{0\}$, let $r_z:=\frac{2r}{|z|}$. Observe that the claim (\ref{eq:main-prop}) of the proposition is trivial if
$|z|\leqslant 2r$. So we may assume that $2r<|z|\leqslant1$ and hence that $r_z\in[2r,1)$. Next let
$n_\ell:=\lceil\frac{\log r_z}{\log\gamma}\rceil$ and $n_{k}:=\lceil\alpha\, n_{k+1}\rceil$ ($k=\ell-1,\dots,0$). 
Observe that $n_0\leqslant n_{1}\leqslant\dots\leqslant n_\ell
\leqslant N:=\lceil\frac{\log(2r)}{\log\gamma}\rceil$. On the other hand,
\begin{displaymath}
n_k\geqslant\alpha n_{k+1}\geqslant\dots\geqslant \alpha^{\ell-k}n_\ell
\geqslant\alpha^{\ell-k}\frac{\log r_z}{\log\gamma}\quad(k=0,\dots,\ell),
\end{displaymath}
so that
\begin{equation}\label{eq:gamma-alpha-r}
\gamma^{n_{k}}=r_z^{\frac{\log\gamma}{\log r_z}n_{k}}
\leqslant r_z^{\alpha^{\ell-k}}\quad\text{and}\quad
\left(\frac{\gamma}{b}\right)^{n_{k-1}}
\leqslant\left(\frac{\gamma}{b}\right)^{\alpha n_{k}}
=\gamma^{n_{k}}\quad\text{for }k=1,\dots,\ell.
\end{equation}
We will also use the fact that $n_{k}\leqslant\alpha n_{k+1}+1$, which yields by induction
\begin{equation*}%\label{eq:other-direction}
n_k\leqslant\alpha^{\ell-k} n_\ell+\frac{1-\alpha^{\ell-k}}{1-\alpha}\quad\text{for }k=\ell,\dots,0\ .
\end{equation*}

For $k=0,\dots,\ell$ define truncated versions of $\Theta_z$,
\begin{displaymath}
\Theta_{z,k}(\xi,x):=
\sum_{n=1}^{n_k}\gamma^n\,s\left(\frac{z}{b^n}\right)
\sin\left(2\pi\left( x_n+\frac{z}{2b^n}\right)\right)\ ,
\end{displaymath}
and, for $z=0$, also rescaled increments
\begin{displaymath}
\begin{split}
\Delta_{z,k}(\xi,x):&=\gamma^{-n_{k-1}}\left(\Theta_{z,k}(\xi,x)-\Theta_{z,k-1}(\xi,x)\right)\\
&=
\sum_{n=1}^{d_k}\gamma^ns\left(\frac{z/b^{n_{k-1}}}{b^n}\right)\sin\left(2\pi\left(x_{n+n_{k-1}}+\frac{z}{2b^{n+n_{k-1}}}\right)\right)
\ ,
\end{split}
\end{displaymath}
where $d_k:=n_k-n_{k-1}$.
Then
\begin{equation}\label{eq:Rz0-est}
\|\Theta_z-\Theta_{z,\ell}\|_\infty
\leqslant
C\gamma^{n_\ell}\ ,
\end{equation}
\begin{equation}\label{eq:Theta0-estimate}
\|\Theta_0\|_\infty\leqslant C\quad\text{and}\quad
\|\Delta_{0,k}-\Theta_0\circ B^{n_{k-1}}\|_\infty
\leqslant C\gamma^{d_k}\ ,
\end{equation}
and
\begin{equation}\label{eq:approx1}
\begin{split}
&\|\Theta_{z,k}-(\Theta_{z,k-1}+\gamma^{n_{k-1}}\Theta_0\circ B^{n_{k-1}})\|_\infty
=
\gamma^{n_{k-1}}\|\Delta_{z,k}-\Theta_0\circ B^{n_{k-1}}\|_\infty\\
\leqslant&
\gamma^{n_{k-1}}\left(\|\Delta_{z,k}-\Delta_{0,k}\|_\infty+C\gamma^{d_k}\right)
\leqslant
C\,\gamma^{n_{k-1}}\left(b^{-n_{k-1}}+\gamma^{d_k}\right)
\leqslant
C\,\gamma^{n_k}
\end{split}
\end{equation}
by (\ref{eq:gamma-alpha-r}).
From now on we denote by $C$ a constant such that these last estimates are  satisfied for $k=1,\dots,\ell$.
Observe that $C$ does neither depend on $r$ nor on $z$.

We have a closer look at measurability properties of $\Delta_{0,k}$, $\Theta_0\circ B^{n_{k-1}}$ and $\Theta_{z,k-1}$:
\begin{enumerate}[$\triangleright$]
\item$\Delta_{0,k}$ depends on $(\xi,x)$ only through
$x_{n_{k-1}+1},\dots,x_{n_k}$, i.e. through $x_{n_{k-1}},k(\xi_{n_{k-1}}),\dots,k(\xi_{n_{k}-1})$.
\item The conditional distribution of $\Theta_0(B^{n_{k-1}}(\xi,x))=-2\pi\sum_{n=1}^\infty\gamma^n\sin(2\pi x_{n+n_{k-1}})$ given $x_{n_{k-1}}$ is $\nu_{x_{n_{k-1}}}$ (Remark~\ref{remark:conditionalTheta}). It has density $h_{x_{n_{k-1}}}$ with respect to Lebesgue measure $m$ (Proposition~\ref{prop:like-Tsujii}).
\item $\Theta_{z,k-1}$ depends on $(\xi,x)$ only through
$x_{1},\dots,x_{n_{k-1}}$, i.e. through $x,k(\xi_{0}),\dots,k(\xi_{n_{k-1}-1})$.
\end{enumerate}

In the rest of the proof we use constants $C_k=(2(\ell-k)+1)C$, $k=0,\dots,\ell$, where $C>0$ is the constant fixed above.
For $j=3,4,5$ let $R_{k,j}=jC+2C_k$ and define probability kernels
\begin{equation*}
\varphi_{k,j}(t)=\frac{1}{2R_{k,j}\gamma^{d_k}}\,1_{[-R_{k,j}\gamma^{d_k},R_{k,j}\gamma^{d_k}]}(t)\ .
\end{equation*}
It is easily checked that for $j=3,4$, $|\eta|\leqslant C\gamma^{d_k}$ and each non-negative $h\in L^1_m(\R)$ holds
\begin{equation}\label{eq:varphi-ineq}
h*\varphi_{k,j}(u+\eta)
\leqslant
%\frac{R_{k,j+1}}{R_{k,j}}
2h*\varphi_{k,j+1}(u)\ .
\end{equation}

Now let
\begin{equation}\label{eq:Erkz-def}
E_{r,k}(z):=\left\{(\xi,x)\in\I^2: h_{x_{n_{k-1}}}*\varphi_{k,4}(\Delta_{0,k}(\xi,x))>r^{-2\delta}
\right\},
\end{equation}
and $G_{r,k}(z):=\I^2\setminus(E_{r,1}(z)\cup\dots\cup E_{r,k}(z))$.
(The argument $z$ reminds of the fact that the sequence $n_0\leqslant n_1 \leqslant\dots\leqslant n_\ell\leqslant N$ depends on $z$ -- but not the common upper bound $N$!)
Observe also that
\begin{enumerate}[$\triangleright$]
\item $1_{E_{r,k}(z)}(\xi,x)$ depends on $(\xi,x)$ only through $x_{n_{k-1}}$ and through 
$\Delta_{0,k}(\xi,x)$, and hence only through $x,k(\xi_0),\dots,k(\xi_{n_k-1})$, so that
\item $1_{G_{r,k-1}(z)}(\xi,x)$ depends on $(\xi,x)$ only through $x,k(\xi_0),\dots,k(\xi_{n_{k-1}-1})$.
\end{enumerate}
Finally
let $E_r$ be the union of all such sets $E_{r,k}(z)$ where the union extends over all $k=0,\dots,\ell$ and all $z\in[2r,1)$. A very crude estimate shows that at most $(\ell+1) N^{\ell+1}$ different sets $E_{r,k}(z)$ can contribute to this union. At the end of the proof we will show that there is a constant $C'>0$ such that
\begin{equation}\label{eq:E-estimate}
m^2(E_{r,k}(z))\leqslant C'\, r^{2\delta}\quad\text{for all such sets }E_{r,k}(z).
\end{equation}
Recalling that  $N=\lceil\frac{\log(2r)}{\log\gamma}\rceil$ and that the choice of $\ell$ did neither depend on $r$ nor on $z$, this shows that $m^2(E_r)\leqslant C''r^{\delta}$ as claimed.

Now let $(\L_x)_{x\in\I}$ be a measurable family of intervals of length $r_z(\approx\gamma^{n_\ell})$, and denote by $\L_x(t)$ the $t$-neighbourhood of $\L_x$. 
Then, by (\ref{eq:Rz0-est}),
\begin{equation}\label{eq:first-step}
m^2\{(\xi,x)\in \I^2\setminus E_r: \Theta_z(\xi,x)\in \L_x\}
\leqslant
m^2\{(\xi,x)\in G_{r,\ell}(z):\Theta_{z,\ell}(\xi,x)\in \L_x(C_\ell\gamma^{n_\ell})\}\ .
\end{equation}
Next we prove for all $k=\ell,\dots,1$ that
\begin{equation}\label{eq:recursive-estimate}
\begin{split}
&m^2\{(\xi,x)\in G_{r,k}(z) : \Theta_{z,k}(\xi,x)\in \L_x(C_k\gamma^{n_k})\}\\
\leqslant& r^{-2\delta}\,4R_{k,3}\gamma^{d_k}\cdot
m^2\{(\xi,x)\in G_{r,k-1}(z) : \Theta_{z,k-1}(\xi,x)\in \L_x(C_{k-1}\gamma^{n_{k-1}})\}\ .
\end{split}
\end{equation}
Indeed, by (\ref{eq:approx1}), (\ref{eq:Theta0-estimate}) and by the measurability properties of $G_{r,k-1}(z)$ and $\Theta_{z,k-1}$,
\begin{equation*}
\begin{split}
&m^2\left\{(\xi,x)\in G_{r,k}(z) : \Theta_{z,k}(\xi,x)\in \L_x(C_k\gamma^{n_k})\right\}\\
\leqslant&
m^2\left\{(\xi,x)\in G_{r,k}(z) : \Theta_{z,k-1}(\xi,x)+\gamma^{n_{k-1}}\Theta_0\circ B^{n_{k-1}}(\xi,x)\in \L_x((C+C_k)\gamma^{n_k})\right\}\\
\leqslant&
m^2\left\{(\xi,x)\in G_{r,k-1}(z)\setminus{E}_{r,k}(z) :
\Theta_{z,k-1}(\xi,x)\in\L_x(C_{k-1}\gamma^{n_{k-1}})\text{ and }\right.\\
&\hspace*{4cm}\left. 
\Theta_0\circ B^{n_{k-1}}(\xi,x)\in \gamma^{-n_{k-1}}\left(\L_x((C+C_k)\gamma^{n_k})-\Theta_{z,k-1}(\xi,x)\right)\right\}\\
=&
\int_{G_{r,k-1}(z)\cap\{\Theta_{z,k-1}\in\L_.(C_{k-1}\gamma^{n_{k-1}})\}}\\
&\hspace*{1cm }
m^2\left(\left\{
\Theta_0\circ B^{n_{k-1}}(\xi,x)\in \J_{(\xi,x)}\right\}
\setminus{E}_{r,k}(z)
\bigg|x,k(\xi_0),\dots,k(\xi_{n_{k-1}-1})\right)
dm^2(\xi,x)\ ,
\end{split}
\end{equation*}
where $\J_{(\xi,x)}=\gamma^{-n_{k-1}}\left(\L_x((C+C_k)\gamma^{n_k})-\Theta_{z,k-1}(\xi,x)\right)$ is an interval of length at most $R_{k,3}\gamma^{d_k}=(3C+2C_k)\gamma^{d_k}$,
and $\tilde\J_x:=\J_{B^{-n_{k-1}}(\xi,x)}$ depends on $(\xi,x)$ only through $x$.
Observing the definition of the set $E_{r,k}(z)$ in (\ref{eq:Erkz-def}), an elementary calculation with conditional probabilities and their behaviour under measure preserving maps shows that the conditional probability in the last integral equals
\begin{equation*}
m^2\left(\left\{
\Theta_0(\xi,x)\in \tilde\J_x\right\}
\cap\left\{h_x*\varphi_{k,4}(\Delta_{0,k}(B^{-n_{k-1}}(\xi,x)))\leqslant r^{-2\delta} \right\}
\bigg|x\right)\circ B^{n_{k-1}}\ .
\end{equation*}
In view of (\ref{eq:Theta0-estimate}) and
(\ref{eq:varphi-ineq}) this can be estimated by
\begin{equation*}
\begin{split}
\leqslant&
m^2\left(\left\{
\Theta_0(\xi,x)\in \tilde\J_x\right\}
\cap\left\{h_x*\varphi_{k,3}(\Theta_0(\xi,x)))\leqslant 2r^{-2\delta} \right\}
\bigg|x\right)\circ B^{n_{k-1}}\\
=&
\int_{\tilde\J_x\cap\{h_x*\varphi_{k,3}\leqslant 2r^{-2\delta}\}}h_x(\theta)\,d\theta\;\bigg|_{x=x_{n_{k-1}}}\\
\leqslant&
\begin{cases}
0&\text{ if }\tilde\J_{x_{n_{k-1}}}\cap\{h_{x_{n_{k-1}}}*\varphi_{k,3}\leqslant 2r^{-2\delta}\}=\emptyset\\
2R_{k,3}\gamma^{d_k}\cdot h_{x_{n_{k-1}}}*\varphi_{k,3}(\tilde{\theta})&
\text{ if }\tilde{\theta}\in\tilde\J_{x_{n_{k-1}}}\cap\{h_{x_{n_{k-1}}}*\varphi_{k,3}\leqslant 2r^{-2\delta}\}
\end{cases}\\
\leqslant&
r^{-2\delta}\,4R_{k,3}\gamma^{d_k}\ ,
\end{split}
\end{equation*}
where we used the fact that the conditional distribution of $\Theta_0$ given $x$ has density $h_x$. This finishes the proof of the recursive estimate (\ref{eq:recursive-estimate}). Combining that estimate with (\ref{eq:first-step}) we obtain
\begin{equation*}
\begin{split}
m^2\{(\xi,x)\in \I^2\setminus E_r: \Theta_z(\xi,x)\in \L_x\}
&\leqslant
\const\,r^{-2\ell\delta}\gamma^{n_\ell-n_0}\\
&\leqslant
\const\,r^{-\eta}\,r_z^{1-\eta}
\leqslant
\const\, r^{1-2\eta}\,|z|^{-(1-\eta)}
\ ,
\end{split}
\end{equation*}
where we used that $\eta=2\ell\delta$ and $n_\ell-n_0\geqslant n_\ell(1-\alpha^\ell)\geqslant n_\ell(1-\eta)$. As the constant does neither depend on $r$ nor on $z$, this proves the claim (\ref{eq:main-prop}) of Proposition~\ref{prop:inductive}.

It remains to prove estimate (\ref{eq:E-estimate}):
\begin{equation*}
\begin{split}
&m^2\left\{(\xi,x)\in\I^2: h_{x_{n_{k-1}}}*\varphi_{k,4}(\Delta_{0,k}(\xi,x))>r^{-2\delta}
\right\}\\
=&
m^2\left\{(\xi,x)\in\I^2: h_{x}*\varphi_{k,4}(\Delta_{0,k}\circ B^{-n_{k-1}}(\xi,x))>r^{-2\delta}
\right\}\\
\leqslant&
r^{2\delta}\int_\I\int_\I h_x*\varphi_{k,4}(\Delta_{0,k}\circ B^{-n_{k-1}}(\xi,x))\,d\xi\,dx\\
\leqslant&
2r^{2\delta}\int_\I\int_\I h_x*\varphi_{k,5}(\Theta_0(\xi,x))\,d\xi\,dx
\quad\text{ by (\ref{eq:Theta0-estimate}) and
(\ref{eq:varphi-ineq})}\\
=&
2r^{2\delta}\int_\I\int_\R h_x*\varphi_{k,5}(\theta)\,h_x(\theta)\,d\theta\,dx\\
\leqslant&
2r^{2\delta}\int_\I \|h_x*\varphi_{k,5}\|_2\,\|h_x\|_2\,dx
\leqslant
2H\,r^{2\delta}\ ,
\end{split}
\end{equation*}
where we used the fact that convolution with a probability kernel cannot increase the $L^2_m$-norm of a function.
This finishes the proof of Proposition~\ref{prop:inductive}.

\end{document}